\documentclass[12pt]{article}
\usepackage{latexsym}
\usepackage{amssymb}
\usepackage{amsmath}

\begin{document}


\newtheorem{theorem}{Theorem}
\newtheorem{problem}{Problem}
\newtheorem{definition}{Definition}
\newtheorem{lemma}{Lemma}
\newtheorem{proposition}{Proposition}
\newtheorem{corollary}{Corollary}
\newtheorem{example}{Example}
\newtheorem{conjecture}{Conjecture}
\newtheorem{algorithm}{Algorithm}
\newtheorem{exercise}{Exercise}
\newtheorem{remarkk}{Remark}

\newcommand{\be}{\begin{equation}}
\newcommand{\ee}{\end{equation}}
\newcommand{\bea}{\begin{eqnarray}}
\newcommand{\eea}{\end{eqnarray}}
\newcommand{\beq}[1]{\begin{equation}\label{#1}}
\newcommand{\eeq}{\end{equation}}
\newcommand{\beqn}[1]{\begin{eqnarray}\label{#1}}
\newcommand{\eeqn}{\end{eqnarray}}
\newcommand{\beaa}{\begin{eqnarray*}}
\newcommand{\eeaa}{\end{eqnarray*}}
\newcommand{\req}[1]{(\ref{#1})}

\newcommand{\lip}{\langle}
\newcommand{\rip}{\rangle}
\newcommand{\uu}{\underline}
\newcommand{\oo}{\overline}
\newcommand{\La}{\Lambda}
\newcommand{\la}{\lambda}
\newcommand{\eps}{\varepsilon}
\newcommand{\om}{\omega}
\newcommand{\Om}{\Omega}
\newcommand{\ga}{\gamma}
\newcommand{\rrr}{{\Bigr)}}
\newcommand{\qqq}{{\Bigl\|}}

\newcommand{\dint}{\displaystyle\int}
\newcommand{\dsum}{\displaystyle\sum}
\newcommand{\dfr}{\displaystyle\frac}
\newcommand{\bige}{\mbox{\Large\it e}}
\newcommand{\integers}{{\Bbb Z}}
\newcommand{\rationals}{{\Bbb Q}}
\newcommand{\reals}{{\rm I\!R}}
\newcommand{\realsd}{\reals^d}
\newcommand{\realsn}{\reals^n}
\newcommand{\NN}{{\rm I\!N}}
\newcommand{\DD}{{\rm I\!D}}
\newcommand{\LL}{{\rm I\!L}}
\newcommand{\degree}{{\scriptscriptstyle \circ }}
\newcommand{\dfn}{\stackrel{\triangle}{=}}
\def\complex{\mathop{\raise .45ex\hbox{${\bf\scriptstyle{|}}$}
     \kern -0.40em {\rm \textstyle{C}}}\nolimits}
\def\hilbert{\mathop{\raise .21ex\hbox{$\bigcirc$}}\kern -1.005em {\rm\textstyle{H}}} 
\newcommand{\RAISE}{{\:\raisebox{.6ex}{$\scriptstyle{>}$}\raisebox{-.3ex}
           {$\scriptstyle{\!\!\!\!\!<}\:$}}} 

\newcommand{\hh}{{\:\raisebox{1.8ex}{$\scriptstyle{\degree}$}\raisebox{.0ex}
           {$\textstyle{\!\!\!\! H}$}}}

\newcommand{\OO}{\won}
\newcommand{\calA}{{\cal A}}
\newcommand{\calB}{{\cal B}}
\newcommand{\calC}{{\cal C}}
\newcommand{\calD}{{\cal D}}
\newcommand{\calE}{{\cal E}}
\newcommand{\calF}{{\cal F}}
\newcommand{\calG}{{\cal G}}
\newcommand{\calH}{{\cal H}}
\newcommand{\calK}{{\cal K}}
\newcommand{\calL}{{\cal L}}
\newcommand{\calM}{{\cal M}}
\newcommand{\calO}{{\cal O}}
\newcommand{\calP}{{\cal P}}
\newcommand{\calX}{{\cal X}}
\newcommand{\calXX}{{\cal X\mbox{\raisebox{.3ex}{$\!\!\!\!\!-$}}}}
\newcommand{\calXXX}{{\cal X\!\!\!\!\!-}}
\newcommand{\gi}{{\raisebox{.0ex}{$\scriptscriptstyle{\cal X}$}
\raisebox{.1ex} {$\scriptstyle{\!\!\!\!-}\:$}}}
\newcommand{\intsim}{\int_0^1\!\!\!\!\!\!\!\!\!\sim}
\newcommand{\intsimt}{\int_0^t\!\!\!\!\!\!\!\!\!\sim}
\newcommand{\pp}{{\partial}}
\newcommand{\al}{{\alpha}}
\newcommand{\sB}{{\cal B}}
\newcommand{\sL}{{\cal L}}
\newcommand{\sF}{{\cal F}}
\newcommand{\sE}{{\cal E}}
\newcommand{\sX}{{\cal X}}
\newcommand{\R}{{\rm I\!R}}
\renewcommand{\L}{{\rm I\!L}}
\newcommand{\vp}{\varphi}
\newcommand{\N}{{\rm I\!N}}
\def\ooo{\lip}
\def\ccc{\rip}
\newcommand{\ot}{\hat\otimes}
\newcommand{\rP}{{\Bbb P}}
\newcommand{\bfcdot}{{\mbox{\boldmath$\cdot$}}}

\renewcommand{\varrho}{{\ell}}
\newcommand{\dett}{{\textstyle{\det_2}}}
\newcommand{\sign}{{\mbox{\rm sign}}}
\newcommand{\TE}{{\rm TE}}
\newcommand{\TA}{{\rm TA}}
\newcommand{\E}{{\rm E\,}}
\newcommand{\won}{{\mbox{\bf 1}}}
\newcommand{\Lebn}{{\rm Leb}_n}
\newcommand{\Prob}{{\rm Prob\,}}
\newcommand{\sinc}{{\rm sinc\,}}
\newcommand{\ctg}{{\rm ctg\,}}
\newcommand{\loc}{{\rm loc}}
\newcommand{\trace}{{\,\,\rm trace\,\,}}
\newcommand{\Dom}{{\rm Dom}}
\newcommand{\ifff}{\mbox{\ if and only if\ }}
\newcommand{\proof}{\noindent {\bf Proof:\ }}
\newcommand{\remark}{\noindent {\bf Remark:\ }}
\newcommand{\remarks}{\noindent {\bf Remarks:\ }}
\newcommand{\note}{\noindent {\bf Note:\ }}

\newcommand{\boldx}{{\bf x}}
\newcommand{\boldX}{{\bf X}}
\newcommand{\boldy}{{\bf y}}
\newcommand{\boldR}{{\bf R}}
\newcommand{\uux}{\uu{x}}
\newcommand{\uuY}{\uu{Y}}

\newcommand{\limn}{\lim_{n \rightarrow \infty}}
\newcommand{\limN}{\lim_{N \rightarrow \infty}}
\newcommand{\limr}{\lim_{r \rightarrow \infty}}
\newcommand{\limd}{\lim_{\delta \rightarrow \infty}}
\newcommand{\limM}{\lim_{M \rightarrow \infty}}
\newcommand{\limsupn}{\limsup_{n \rightarrow \infty}}

\newcommand{\ra}{ \rightarrow }

\newcommand{\ARROW}[1]
  {\begin{array}[t]{c}  \longrightarrow \\[-0.2cm] \textstyle{#1} \end{array} }

\newcommand{\AR}
 {\begin{array}[t]{c}
  \longrightarrow \\[-0.3cm]
  \scriptstyle {n\rightarrow \infty}
  \end{array}}

\newcommand{\pile}[2]
  {\left( \begin{array}{c}  {#1}\\[-0.2cm] {#2} \end{array} \right) }

\newcommand{\floor}[1]{\left\lfloor #1 \right\rfloor}

\newcommand{\mmbox}[1]{\mbox{\scriptsize{#1}}}

\newcommand{\ffrac}[2]
  {\left( \frac{#1}{#2} \right)}

\newcommand{\one}{\frac{1}{n}\:}
\newcommand{\half}{\frac{1}{2}\:}

\def\le{\leq}
\def\ge{\geq}
\def\lt{<}
\def\gt{>}

\def\squarebox#1{\hbox to #1{\hfill\vbox to #1{\vfill}}}
\newcommand{\qed}{\hspace*{\fill}
           \vbox{\hrule\hbox{\vrule\squarebox{.667em}\vrule}\hrule}\bigskip}

\title{Some  remarks about the positivity of  random variables on a Gaussian
  probability space }

\author{D. Feyel and A. S. \"Ust\"unel}
\date{ }
\maketitle
\noindent
{\bf Abstract:}{\small{ Let $(W,H,\mu)$ be an abstract Wiener space
    and let $L\in \LL\log \LL(\mu)$ is a positive random
    variable. Using the measure transportation  of  Monge-Kantorovitch, we
    prove that operator corresponding to  the kernel of the projection
    of $L$ on the second Wiener chaos is lower bounded by a
    semi-positive Hilbert-Schmidt operator. }}\\

\vspace{0.5cm}
\noindent
{{\bf Quelques  remarques  sur la positivit\'e des variables
    al\'eatoires d\'efinies sur un espace gaussien}} \\ 

\vspace{0.5cm}
\noindent
{\bf Resum\'e:}{\small{ Soit $(W,H,\mu)$ un espace de Wiener abstrait
    et soit $L\in\LL\log\LL$ une variable al\'eatoire positive. A
    l'aide de la th\'eorie de transport de mesure de
    Monge-Kantorovitch, nous montrons que le noyau de la projection de
    $L$ dans le second chaos de Wiener est un op\'erateur de spectre 
    inf\'erieurement born\'e et l'op\'erateur correspondant est
    inf\'erieurement born\'e par un op\'erateur Hilbert-Schmidt semi-positif.
}}

\section{Version fran\c{c}aise abr\'eg\'ee}
Soit $(W,H,\mu)$ un espace de Wiener abstrait: W est un Fr\'echet s\'eparable
localement convexe, $\mu$ est une  mesure gaussienne dont le support
set $W$ et $H$
est l'espace de Cameron-Martin dont le produit scalaire et la norme sont
not\'es respectivement  $(\cdot,\cdot)_H$ et $|\cdot|_H$. On notera
par $\nabla$ la fermeture par rapport \`a $\mu$ de  la d\'eriv\'ee
 dans la direction de $H$. En particulier, pour un espace 
hilbertien $M$, $\DD_{2,k}(M)$ est  l'espace de classes
d'\'equivalences de  fonctions mesurables, \`a valeurs dans $M$, dont
les d\'eriv\'ees d'ordre $k\in \NN$ sont de carr\'e int\'egrables par
rapport \`a la norme du produit tensoriel Hilbert-Schmidt $M\otimes
H^{\otimes k}$, o\`u $H^{\otimes k}$ est l'espace des $k$-tenseurs
Hilbert-Schmidt; si $M=\R$ alors nous noterons $\DD_{2,k}$ au lieu de
$\DD_{2,k}(\R)$ (cf.\cite{F-P}, \cite{Mall},  \cite{ASU}).  
On notera par $\delta$ l'adjoint de $\nabla$ par rapport \`a $\mu$,
qui est une application continue de $\DD_{2,1}(M\otimes H^{\otimes
  k+1})$ dans $\DD_{2,1}(M\otimes H^{\otimes k})$.  Noter  que
$\delta\circ \nabla$ est  l'op\'erateur d'Ornstein-Uhlenbeck,
il sera not\'e  $\calL$. A l'aide de l'in\'egalit\'e de Meyer, on peut
d\'efinir les espaces de Sobolev d'ordre n\'egatif
$(\DD_{p,\alpha},\,\alpha\in \R, p>1)$ et on note
$\DD'=\cup_{p>1,\alpha\in \R}\DD_{p,\alpha}$, qui est dual de l'espace
$\DD=\cap_{p>1,\alpha\in \R}\DD_{p,\alpha}$ (cf.\cite{Mall,ASU}).

Quand $W$ est l'espace de Wiener classique, i.e., $W=C_0([0,1],\R)$,
$H=H_1([0,1],dt)$ (i.e., les primitives des \'el\'ements de  $\LL^2([0,1],dt)$)
il est bien connu que chaque \'el\'ement $L$ de 
$\LL^2(\mu)$ admet une d\'ecomposition unique comme 
$$
L=E[L]+\sum_{n=1}^\infty I_n(L_n)\,,
$$
o\`u $L_n\in \LL_s^{2}([0,1]^n)$ et ce dernier repr\'esente les
fonctions sym\'etriques et de carr\'e int\'egrables sur
$[0,1]^n$. Soit $H^{\odot n}$ repr\'esente le produit tensoriel
sym\'etrique d'ordre $n$ de $H$, qui est isomorphe \`a
$\LL_s^{2}([0,1]^n)$. Si on note par $i_n,\,n\geq 1$ cet isomorphisme,
on peut montrer facilement que $I_n(L_n)=\delta^n(i_n(L_n))$, o\`u
$\delta^n=(\nabla^n)^\star$ par rapport \`a $\mu$. Avec ces relations,
on peut montrer \`a partir de la formule de Taylor que 
$$
L=E[L]+\sum_{n=1}^\infty \frac{1}{n!}\delta^n (E[\nabla^n L])\,,
$$
cf. \cite{Mc-Kean},  \cite{Str} et aussi \cite{ASU,BOOK}.

Soit $\nu$ une autre probabilit\'e, notons par  $\Sigma(\mu,\nu)$
l'ensemble des probabilit\'es sur $W\times W$ de  marginales
$\mu$ et $\nu$. On note $J$ la fonctionnelle d\'efinie sur
$\Sigma(\mu,\nu)$ par  $J(\beta)=\int_{W\times
  W}|x-y|_H^2d\beta(x,y)$. Dans le cas o\`u $W$ est de dimension
finie, le probl\`eme  de Monge-Kantorovitch consiste \`a  trouver
une mesure $\ga\in \Sigma(\mu,\nu)$ telle que la distance de Wasserstein
$$
d^2_H(\mu,\nu)=\inf\{J(\beta):\,\beta\in \Sigma(\mu,\nu)\}
$$
soit atteinte en $\ga$. Ce probl\`eme a \'et\'e r\'esolu dans
\cite{BRE} en dimension finie (cf. aussi \cite{Feyel} pour un survol
rapide). Nous l'avons r\'esolu dans \cite{F-U1,F-U2} (c.f. aussi
\cite{F-U3})  quand la  dimension de $H$ est
infinie. Expliquons plus pr\'ecis\'ement le cas particulier qui sera
utilis\'e dans cette note: si $\nu$ est de la forme $d\nu=Ld\mu$,
alors il existe une fonction $\varphi$, appel\'ee le potentiel de
transport,  appartenant \`a  $\DD_{2,1}$, telle que
$T:W\to W$ d\'efinie par $T=I_W+\nabla\varphi$ satisfasse  $T\mu=\nu$
et telle que  $\ga=(I_W\times T)\mu$ soit  l'unique mesure dans
$\Sigma(\mu,\nu)$ 
satisfaisant  $J(\ga)=d_H^2(\mu,\nu)$. De plus $\varphi$ est $1$-convexe: 
une variable al\'eatoire $f:W\to \R\cup\{\infty\}$ est dite
$r$-convexe, $r\in \R$,  si $h\to
\frac{r}{2}|h|_H^2+f(w+h)$ est convexe sur $H$ \`a valeurs dans
$\L^0(\mu)$ (\cite{F-U}); si $r=0$, on l'appelle $H$-convexe. De
m\^eme $f$ s'appelle $H$-concave ou $H$-log-concave  si,
respectivement  $-f$ est $H$-convexe ou $-log f$ est $H$-convexe.   
Avec les hypoth\`eses ci-dessus $T$ admet un inverse p.s.,  
 not\'e $S$,  de la forme $S=I_W+\eta$. De plus si $\nabla$
est fermable par rapport \`a $\nu$ alors $\eta:W\to H$ est de la forme
$\eta=\nabla\psi$ o\`u  $\psi\in L^2(\nu)$ est $\nu$-differentiable dans la
direction de $H$.
Notons que nous avons d\'ej\`a d\'emontr\'e dans \cite{F-U4} que 
 $\varphi$ est un \'el\'ement de $\DD_{2,2}$ au
lieu de $\DD_{2,1}$ si la densit\'e $L\in \LL\log\LL$ est
$H$-log-concave. Ce qui rend possible le calcul du  jacobien 
$$
\La=\dett(I_H+\nabla^2\varphi)
\exp\left\{-\calL\varphi-\frac{1}{2}|\nabla\varphi|_H^2\right\}\,,
$$
o\`u $\dett(I_H+\nabla^2\varphi)$ est  le d\'eterminant
modifi\'e de Carleman-Fredholm (cf.\cite{Du-S,BOOK}).

\section{Main results}
Here is the first notable  result of this note:
\begin{theorem}
\label{+-thm}
Assume that $L\in \LL^2(\mu)$ is a positive random variable and let
$\varphi$ be the forward potential function associated to the
Monge-Kantorovitch problem in $\Sigma(\mu,\nu)$, where
$d\nu=\frac{1}{E[L]}Ld\mu$.  Then 
the following operator inequality holds true:
\begin{equation}
\label{ineq-1}
\frac{1}{2E[L]}\left\{E[\nabla^2L]-\frac{E[\nabla L]\otimes E[\nabla
    L]}{E[L]}\right\}\geq E[\nabla^2\varphi]\,. 
\end{equation}
\end{theorem}
\proof
Let us note first that, even if $\varphi$ is not in $\DD_{2,2}$, the
term $E[\nabla^2\varphi]$ is a well-defined Hilbert-Schmidt operator
since the constants are the elements of the space of the test
functions $\DD=\cap_{p,k}\DD_{p,k}$. 
Without loss of generality, we may assume that $E[L]=1$. Let then
$\nu$ be the measure $d\nu=Ld\mu$. Since $E[L\log L]<\infty$, the
Wasserstein distance $d_H(\mu,\nu)<\infty$, consequently, there exists
a $1$-convex map  $\varphi\in \DD_{2,1}$ such that the transformation
$T=I_W+\nabla\varphi$ solves the problem of Monge and the measure
$(I\times T)\mu$ is the unique solution of Monge-Kantorovitch problem
on $\Sigma(\mu,\nu)$. For an $h\in H$, let $\rho(\delta h)$ denote the
Wick exponential $\rho(\delta h)=\exp(\delta
h-\frac{1}{2}|h|_H^2)$. For any $t\in \R$, we have 
\beaa
E[L\,\rho(\delta(th))]&=&E[\rho(\delta(th))\circ T]\\
&=&E\left[(\exp\left(t\delta h-\frac{t^2}{2}|h|_H^2\right)\circ
  T\right]\\
&=&E\left[\exp\left(t\delta
    h+t(\nabla\varphi,h)_H-\frac{t^2}{2}|h|_H^2\right)\right]\,.
\eeaa
A first order  differentiation of this equality at $t=0$ gives that 
$$
E[(\nabla L,h)_H]=E[(\nabla\varphi,h)_H]\,,
$$
for any $h\in H$, hence
\begin{equation}
\label{1-st-order}
E[\nabla L]=E[\nabla\varphi]\,.
\end{equation}
The second order differentiation at $t=0$ and the integration by parts
formula, which follows from the fact that $\delta=\nabla^\star$, gives 
\beaa
E\left[(\nabla^2L,h\otimes h)_2\right]&=&E\left[(\delta
h+(\nabla\varphi,h)_H)^2-|h|_H^2\right]\\
&=&E\left[2\,\delta h\,(\nabla\varphi,h)_H+(\nabla\varphi,h)_H^2\right]\\
&=&E\left[2(\nabla^2\varphi,h\otimes h)_2+(\nabla\varphi,h)_H^2\right]\,,
\eeaa
for any $h\in H$, where $(\cdot,\cdot)_2$ denotes the Hilbert-Schmidt
scalar product. Hence
 combining this  with the relation (\ref{1-st-order}) gives 
\beaa
E[\nabla^2L]&=&E[\nabla\varphi\otimes \nabla\varphi]+2E[\nabla^2\varphi]\\
&\geq&E[\nabla\varphi]\otimes E[\nabla\varphi]+2E[\nabla^2\varphi]\\
&=&E[\nabla L]\otimes E[\nabla L]+2E[\nabla^2\varphi]
\eeaa
 and  the inequality (\ref{ineq-1}) follows.
\qed

\begin{remark}
Note that the inequality of Theorem \ref{+-thm} is different in spirit than the
results of \cite{R-M}.
\end{remark}

\noindent
We can extend the inequality (\ref{ineq-1}) as follows: 
\begin{corollary}
\label{cor-1}
Assume that $m\in \DD'$ is a positive distribution and denote again by
$m$ the Radon measure on $W$ which corresponds to it
(cf. \cite{ASU}). Let $m_2$ be the projection of the distribution $m$
to the second Wiener chaos, which is equal to
$\frac{1}{2}\delta^2M_2$, where $M_2$ is the element
of $H\otimes H$ defined by $M_2(h\otimes k)=\langle
m,\delta^2(h\otimes k)\rangle$.  If the Wasserstein distance
$d_H(\mu,m)$ is finite, we have again 
\begin{equation}
\label{ineq-3}
\frac{1}{2m(W)}\left\{M_2-\frac{M_1\otimes M_1}{m(W)}\right\}\geq
E[\nabla^2\varphi]\,, 
\end{equation}
where $M_1\in H$ is defined by $(M_1,h)_H=\langle m,\delta h\rangle$,
$h\in H$. 
\end{corollary}
\proof
It suffices to apply Theorem \ref{+-thm} to the case $P_tm$, where
$P_t$ is the Ornstein-Uhlenbeck semigroup. Then, from \cite{F-U3}, the
corresponding transport map $\varphi_t$ converges to the transport map
$\varphi$ corresponding to the Monge-Kantorovitch problem for
$\Sigma(\mu,m)$ in $\DD_{2,1}$,  as $t\to 0$.
\qed

\noindent
We have also a weaker inequality whose difference with respect to
(\ref{ineq-1}) is that the Hilbert-Schmidt operator
$2E[\nabla^2\varphi]$ is replaced by the identity operator of $H$:
\begin{proposition}
\label{feyel-prop}
For any positive random variable $L\in \LL^2(\mu)$,  the following
inequality is valid:
$$
I_H+\frac{1}{E[L]}E[\nabla^2L]\geq \frac{1}{E[L]^2}E[\nabla L]\otimes
E[\nabla L]\,,
$$
where $I_H$ denotes the identity operator of $H$. In particular the
projection of $L$ in the second order Wiener chaos divided by the
expectation of $L$ is $1$-convex.
\end{proposition}
\proof
Again, we may suppose that $E[L]=1$. Let $l(t)=E[L\,\rho(\delta(th))$,
$h\in H$, then we have 
\begin{equation}
\label{2nd-deg}
l'(0)^2\leq |h|_H^2+l''(0)
\end{equation}
To see this
inequality it suffices to remark that $\la^2-2\la\delta h+(\delta
h)^2\geq 0$ for any $\la\geq 0$, hence taking the expectation with
respect to $Ld\mu$ is again positive. Hence the discriminant of the
second order polynomial in $\la$ should be negative and this proves
(\ref{2nd-deg}). To complete the proof of the proposition, it suffices
to remark that $l'(0)=E[\nabla_hL]$ and that
$l''(0)=\trace\left(E[\nabla^2L]\,(h\otimes h)\right)$. The
$1$-convexity of $\delta^2\left(\frac{E[\nabla^2L]}{2}\right)$ is immediate.
\qed

\noindent
Proposition \ref{feyel-prop} extends also to the positive elements of
$\DD'$ and it is to be noted that in this case we do not need the
hypothesis about the finiteness of the Wasserstein distance:
\begin{corollary}
Assume that $m\in \DD'$ is a positive distribution and denote again by
$m$ the Radon measure on $W$ which corresponds to it
(cf. \cite{ASU}). Using  the notations of Corlollary \ref{cor-1}, we
have again  
\begin{equation}
\label{ineq-5}
I_H+\frac{1}{m(W)}M_2\geq \frac{1}{m(W)^2}M_1\otimes M_1\,.
\end{equation}
In particular, the projection of $m$ in the second Wiener chaos is
$1$-convex. 
\end{corollary}
\proof
It suffices to regularize again $m$ with $P_t$ and then pass to the
limit as $t\to 0$.
\qed

\vspace{2cm}

{\footnotesize}
\begin{tabular}{ll}
D. Feyel & A.S. \"Ust\"unel,\\
Universit\'e d'Evry-Val-d'Essone & ENST, D\'ept. Infres\\
D\'ept. de Math\'ematiques & 46 Rue Barrault,,  \\
91025 Evry Cedex& 75013 Paris,  \\
France&France\\
feyel@maths.univ-evry.fr&ustunel@enst.fr 
\end{tabular}

\end{document}